\documentclass{elsarticle}
\usepackage[utf8]{inputenc}
\usepackage{amssymb}
\usepackage{csquotes}
\usepackage{amsmath}
\usepackage{hyperref}
\usepackage{amsthm}
\usepackage{color}
\newtheorem{theorem}{Theorem}

\newcommand{\coloneqq}{\mathrel{\mathop:}=}
\renewcommand{\vec}[1]{\mathbf{#1}}

\def\aaa#1{{\color{black}#1}}

\journal{Systems and Control Letters}

\begin{document}
\begin{frontmatter}
\title{A Globally Asymptotically Stable Polynomial Vector Field with Rational Coefficients and no Local Polynomial Lyapunov Function}
 \author{Amir Ali Ahmadi and Bachir El Khadir \corref{fundinginfo}}
\cortext[fundinginfo]{ The authors are with  the department of Operations Research and Financial Engineering at Princeton University. Email:
  \{\texttt{a\_a\_a}, \texttt{bkahdir}\}\texttt{@princeton.edu}. This work was partially supported by the DARPA Young Faculty Award, the Young Investigator Award of the AFOSR, the CAREER Award of the NSF, the Google Faculty Award, and the Sloan Fellowship.
}
\address{Princeton University}
\begin{abstract}
We give an explicit example of a two-dimensional polynomial vector field of degree seven that has rational coefficients, is globally asymptotically stable, but does not admit an analytic Lyapunov function even locally.
\end{abstract}
\end{frontmatter}

\section{Introduction and Motivation}
\label{sec:introduction}

We are concerned in this paper with a continuous time dynamical system
\begin{equation}
  \label{eq:autonomous.dynamical.system}
  \dot {\vec x} = f({\vec x}),
\end{equation}
where $f: \mathbb R^n \rightarrow \mathbb R^n$ is a \emph{polynomial} and has an equilibrium point at the origin, i.e., $f(0) = 0$. Polynomial differential equations appear throughout engineering and the sciences and the study of stability of their equilibrium points has been a problem of long-standing interest to mathematicians and control theorists.

We recall that the origin of \eqref{eq:autonomous.dynamical.system} is said to be a \emph{locally asymptotically stable} (LAS) equilibrium if it is \emph{stable in the sense of Lyapunov} (i.e., if for every \(\epsilon>0\), there exists a \(\delta=\delta(\epsilon)>0\) such that
$\|\vec x(0)\|<\delta\ \Rightarrow \|\vec x(t)\|<\epsilon$ \aaa{for all} $t\geq0$) and if there exists a scalar \(\hat \delta > 0\) such that $$\|\vec x(0)\|< \hat \delta \ \Rightarrow \ \lim_{t\rightarrow\infty}\vec x(t)=0.$$
We say that the origin of \eqref{eq:autonomous.dynamical.system} is a \emph{globally asymptotically stable} (GAS) if it is stable in the sense of Lyapunov and if \(\lim_{t\rightarrow\infty}\vec x(t)=0\) for any initial condition $x(0)$ in \(\mathbb{R}^n\).

We also recall \aaa{(see, e.g.,~\cite{khalil_nonlinear_2001})} that the origin of \eqref{eq:autonomous.dynamical.system} is LAS if there exists a continuously differentiable (Lyapunov) function $V:\mathbb R^n \rightarrow \mathbb R$ that vanishes at the origin and satisfies $V({\vec x}) > 0$ and $-\langle \nabla V({\vec x}), f({\vec x}) \rangle > 0$ for all ${\vec x} \in S \setminus \{0\}$, where $S$ is a neighborhood of the origin. Moreover, if $V$ is in addition radially unbounded (i.e., satisfies  $V(\vec x) \rightarrow \infty$ when $\|\vec x\| \rightarrow \infty$) and if $S = \mathbb R^n$, then the origin is GAS. We call a function satisfying the former (resp. the latter) requirements a \emph{local (resp. global) Lyapunov function}. \aaa{It is also well known that existence of such Lyapunov functions is not only sufficient, but also necessary for local/global asymptotic stability~\cite{khalil_nonlinear_2001}.}

Since the vector field in \eqref{eq:autonomous.dynamical.system} is polynomial, it is natural to search for Lyaponuv functions that are polynomials themselves. This approach has become widely popular in the last couple of decades due to \aaa{the advent of optimization-based algorithms that automate the search for a polynomial Lyapunov function. Arguably, the most prominent such algorithm is based on \emph{sum of squares optimization,} which reduces this search to a semidefinite program~\cite{parrilo_structured_2000,papachristodoulou_construction_2002,henrion_positive_2005,jarvis-wloszek_controls_2003,chesi_lmi_2010,henrion_guest_2009,zheng2017exploiting}. Alternatives to this approach that are based on linear programming or other algebraic techniques have also appeared in recent years~\cite{kamyar2014polynomial,dsosArxiv,kamyar2013solving,ben2015stability}.} 
\aaa{As the algorithmic construction of polynomial Lyapunov functions has been the focus of intense research in recent years, it is natural to ask whether existence of a Lyapunov function within this class is guaranteed. This is the case, e.g., if the goal is to prove exponential stability of an equilibrium point over a bounded region~\cite{peet_exponentially_2009},~\cite{peet2012converse}. Our focus in this paper, however, is on the basic question of whether 
asymptotic stability of an equilibrium point implies existence of a polynomial Lyapunov function.} As is well known, the answer is positive when the degree of the vector field in \eqref{eq:autonomous.dynamical.system} is equal to one. Indeed, asymptotically stable linear systems always admit a quadratic Lyapunov function.




Unlike the linear case, stable polynomial vector fields of degree as low as 2 may fail to admit a polynomial Lyapunov function. Indeed, in \cite{ahmadi_globally_2011}, it is shown that the simple vector field 

\begin{equation}
  \label{eq:dynamical-system-no-rational-lyaponuv}
  \begin{array}{lll}
    \dot{x}&=&-x+xy \\
    \dot{y}&=&\ \ -y
  \end{array}
\end{equation}

\noindent is globally asymptotically stable (e.g. as certified by the Lyapunov function $V(x, y) = \log(1+x^2) + y^2$),
but does not admit a (global) polynomial Lyapunov function. Note however, that the linearization of \eqref{eq:dynamical-system-no-rational-lyaponuv} around the origin is asymptotically stable, and hence this nonlinear system admits a local quadratic Lyapunov function.

In \cite[Prop. 5.2]{bacciotti_liapunov_2006}, Bacciotti and Rosier show that the vector field

\begin{equation}
  \label{eq:Bacciotti.Rosier.f0}
  \begin{aligned}
  \left(\begin{array}{l}
          \dot x\\\dot y
        \end{array}\right) =
  &\left(\begin{array}{lll}
          -2\lambda y(x^2+y^2)-2y(2x^2+y^2) \\
          4\lambda x(x^2+y^2)+2x(2x^2+y^2)
        \end{array}\right)
      \\&- (x^2+y^2)
      \left(\begin{array}{lll}
                         4\lambda x(x^2+y^2)+2x(2x^2+y^2)\\
                         2\lambda y(x^2+y^2)-2y(2x^2+y^2)
            \end{array}\right)
        \end{aligned}
      \end{equation}
is globally asymptotically stable for any scalar $\lambda \ge 0$  (e.g. as certified by the Lyapunov function $V_{\lambda}(x, y) = (x^2+y^2)(2x^2+y^2)^{\lambda}$), but does not admit a \emph{local} polynomial Lyapunov function for any $\lambda$ which is \emph{irrational}.\footnote{In fact, they show that for irrational $\lambda$, the system \eqref{eq:Bacciotti.Rosier.f0} does not even admit a local analytic Lyaponuv function.} However, the validity of this statement crucially relies on the parameter $\lambda$ being irrational. Indeed, for any rational value of $\lambda\geq 0$, the system admits a global polynomial Lyaponuv function, which is e.g. simply an appropriate integer power of $V_{\lambda}$.


Our contribution in this paper is to give an example of a (globally) asymptotically stable polynomial vector field with \emph{rational coefficients} that does not admit a \emph{local} polynomial (or even analytic) Lyaponuv function. Our construction is inspired by and is similar to that of Bacciotti and Rosier~\cite{bacciotti_liapunov_2006}. However, by adapting their underlying proof technique, we are able to prove stability with a Lyapunov function which is the ratio of two polynomials. This allows us to use only rational coefficients in the construction of the vector \aaa{field.\footnote{Note that by rescaling, one can always change a polynomial vector field with rational coefficients to a polynomial vector field with integer coefficients without changing the properties of stability or validity of a candidate Lyapunov function.}}




Our interest in studying polynomial vector fields with rational coefficients partly stems from the fact that in practice, most (if not all) vector fields that are analyzed on a computer (e.g. by an optimization-based algorithm) have rational coefficients. Therefore, if it was true that such vector fields always had polynomial Lyapunov functions, one could restrict attention to this function class for all practical purposes and use techniques such as sum of squares optimization to algorithmically find these Lyapunov functions. Because of this practical motivation, existence of the counterexample that we present in this paper was regarded as a significant unresolved question in the community; see e.g. the ending paragraph in~\cite[Sect. IV]{rational_stability}.


Polynomial vector fields with rational coefficients are also important from the viewpoint of complexity analysis in the standard Turing model. For example, it is not known whether the problem of testing local asymptotic stability is decidable for this class of vector fields. Indeed, this is an outstanding open problem suggested by Arnold, which appears e.g. in \cite{doria_arnolds_1993}, \cite{arnold_problems_1976}:
\begin{displayquote}
``Let a vector field be given by polynomials of a fixed degree, with rational coefficients. Does an algorithm exist, allowing to decide, whether the stationary point is stable?''
\end{displayquote}
In \cite{doria_arnolds_1993}, Arnold is quoted to have conjectured that the answer to the above question is negative:
\begin{displayquote}
``My conjecture has always been that there is no algorithm for some sufficiently high degree and dimension.''
\end{displayquote}

This conjecture also motivates the example in this paper: if it was true that LAS polynomial vector fields with rational coefficients always admitted polynomial Lyaponuv functions of a computable degree, then 
the problem of testing stability would become decidable. This is because one can e.g. use the quantifier elimination theory of Tarski and Seidenberg \cite{tarski_decision_1998},~\cite{seidenberg_new_1954} to test, in finite time, whether a polynomial vector field admits a local polynomial Lyaponuv function of a given degree.

We end our introduction by noting that, interestingly, there is a parallel to these questions in the study of switched linear systems in discrete time. There, the problem of testing asymptotic stability is similarly not known to be decidable
\cite[Problem 10.2]{blondel_unsolved_2004},~\cite{jungers_joint_2009}. One can show, however, that if the so called ``finiteness conjecture'' \cite{lagarias_finiteness_1995} is true for rational matrices, then asymptotic stability becomes decidable. This conjecture is known to be false over the reals \cite{hare_explicit_2011}, but is currently unresolved for rational matrices \cite{jungers_finiteness_2008}.

\section{The Main Result}         
Our contribution in this paper is to prove the following theorem.
\begin{theorem}
  \label{thm:poly.vector.field.with.rat.coeffs.but.no.poly.lyap}
  The polynomial vector field
  \begin{equation}
    \label{eq:poly.vector.field.with.rat.coeffs}
    \begin{pmatrix}\dot x\\\dot y\end{pmatrix} = f(x, y),
  \end{equation}
  with
  \begin{equation*}
      f(x, y) =  \begin{pmatrix}- 2y(- x^4 + 2x^2y^2 + y^4)\\2x(x^4 + 2x^2y^2 - y^4)\end{pmatrix} - (x^2+y^2)\begin{pmatrix}2x(x^4 + 2x^2y^2 - y^4)\\ 2y(- x^4 + 2x^2y^2 + y^4)\end{pmatrix},
    \end{equation*}
  is globally asymptotically stable but does not admit an analytic Lyapunov function even locally.
\end{theorem}

\vspace{-5mm}

\begin{figure}[h]
\centering
\includegraphics[width=0.5\textwidth]{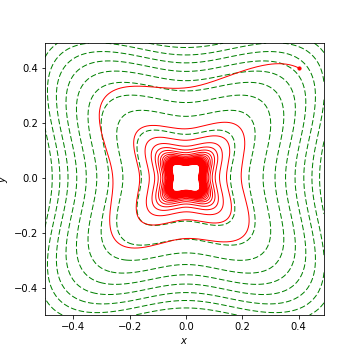}
\caption{A typical trajectory of the vector field in \eqref{eq:poly.vector.field.with.rat.coeffs} and the level sets of the Lyapunov function W in \eqref{eq:rational.lyap}.}
\end{figure}

\begin{proof}
We prove that the vector field in (\ref{eq:poly.vector.field.with.rat.coeffs}) is globally asymptotically stable by means of the rational Lyapunov function \aaa{defined as   \begin{equation}
    \label{eq:rational.lyap}
    W(x, y) = \frac{x^4+y^4}{x^2+y^2}~ \forall (x,y) \neq (0,0), \text{ and } W(0,0)=0.
  \end{equation}}
\aaa{Note} that the function $W$ is continuously differentiable \aaa{on $\mathbb R^2$}, positive definite \aaa{(i.e., satisfies $W(x,y)>0$ for all $(x,y)\neq (0,0)$)}, and radially unbounded. \aaa{Radial unboundedness can be seen, e.g., by noting that since $||(x,y)^T||_2 \leq 2^{1/4} ||(x,y)^T||_4$ for all $(x,y) \in \mathbb{R}^2$, we have $$W(x,y)=\frac{||(x,y)^T||_4^4}{||(x,y)^T||_2^2} \geq \frac12 ||(x,y)^T||_2^2,~\forall (x,y) \in \mathbb{R}^2.$$}
Let us examine the gradient of $W$. A straightforward calculation gives $$\nabla W(x, y) = \frac1{(x^2+y^2)^2} \begin{pmatrix}a(x, y)\\b(x, y)\end{pmatrix},$$
where
\(a(x, y) = 2x(x^4 + 2x^2y^2 - y^4)\)
and
\(b(x, y) = 2y(- x^4 + 2x^2y^2 + y^4)\).

If we let \(f_0 = \begin{pmatrix}-b\\a\end{pmatrix}\), and \(f_1 = -(x^2+y^2) \begin{pmatrix}a\\b\end{pmatrix}\), then $f = f_0 + f_1$, and
\begin{align*}
  \langle \nabla W, f \rangle &= \langle \nabla W, f_0 \rangle + \langle \nabla W, f_1 \rangle
  \\&= 0 - \frac{a^2+b^2}{x^2+y^2}.
\end{align*}
\aaa{We show that $\langle \nabla W, f \rangle$ is negative when $(x,y)\neq (0,0)$ by observing} that for every $(x, y) \in \mathbb R^2 \setminus \{(0, 0)\}$, $a(x, y)$ and $b(x, y)$ cannot both be zero. Indeed, if $a(x, y) = b(x, y) = 0$ for some $(x, y) \in \mathbb R^2$, then $$ya(x, y) + xb(x, y) = 8(xy)^3 = 0,$$ therefore $x = 0$ or $y = 0$. If $x = 0$ for example (the case $y = 0$ is similar), then $b(x, y) = 2 y^5$, and hence $y = 0$ as well. This shows that 
$$\langle \nabla W(x, y), f(x, y) \rangle < 0 \quad \forall (x, y) \ne (0, 0),$$
and hence $W$ is a global Lyapunov function which proves that the vector field is GAS.

Let us now show that $f$ does not admit an analytic Lyapunov function locally. Assume for the sake of contradiction that such a function $p: \mathbb R^2 \rightarrow \mathbb R$ exists. By analyticity, \(p = \sum_{k =0}^\infty p_k\), where $p_k$ is a homogeneous polynomial of degree $k$. Let $p_{k_0}$ be the first non-vanishing term. \aaa{Note that $k_0 \ge 2$ as $$p(0,0)=0\Rightarrow p_0=0,$$  $$p\geq 0, p(0,0)=0\Rightarrow \nabla p(0,0)=(0,0)^T\Rightarrow p_1(x,y)=0,\forall (x,y)\in\mathbb{R}^2.$$
Here, the first implication follows from the fact that the origin is a global minimum for $p$. Observe now that}
\begin{align*}
\langle  \nabla p, f \rangle &=  \langle  \nabla  \sum_{k =\aaa{k_0}}^\infty  p_k , f_0 +f_1\rangle
\\&=  \langle \nabla p_{k_0}, f_0  \rangle + q,
\end{align*}
\noindent where \aaa{$q \mathrel{\mathop{:}}=\langle \nabla p_{k_0},f_1 \rangle +\sum_{k=k_0+1}^\infty \langle \nabla p_k, f_0+f_1 \rangle$. Note that all terms in $q$ have degree higher than the degree of the (homogeneous) polynomial $\langle \nabla p_{k_0}, f_0\rangle$. This is because $f_1$ has higher degree than $f_0$ and the index of the sum in the definition of $q$ starts at $k_0+1$.} \aaa{Since $\langle \nabla p,f \rangle \le 0$ (as we are assuming that $p$ is a Lyapunov function), and since  $\langle \nabla p_{k_0}, f_0\rangle$ constitutes the terms of $\langle \nabla p,f \rangle $ of lowest order, it must be that $\langle \nabla p_{k_0}, f_0\rangle$ is nonpositive in a small enough neighborhood of the origin. But as $\langle \nabla p_{k_0}, f_0\rangle$ is homogeneous, this implies that} 
\begin{equation}\label{eq:pk0dot<=0}
 \langle \nabla p_{k_0}(x, y), f_0(x, y)  \rangle \le 0 \quad \forall (x, y) \in \mathbb R^2.
\end{equation}


We now claim that \aaa{the (homogeneous) polynomial} $p_{k_0}$ must be constant on the 1-level set of $W$, which we denote by
$$M \coloneqq \{(x, y) \in \mathbb R^2 \; | \; W(x, y) = 1\}.$$

\aaa{Since $W$ is continuous (resp. radially unbounded), it follows that $M$ is closed (resp. bounded).} In addition, $f_0$ is continuously differentiable and does not vanish on $M$, \aaa{as we have already argued that $a(x,y)$ and $b(x,y)$ cannot simultaneously vanish except at the origin}. Moreover, trajectories of the vector field $f_0$ that start in $M$ remain in $M$ as \aaa{one can verify that}

$$\langle \nabla W, f_{0} \rangle = 0.$$
Hence, by the Poincar\'e-Bendixson Criterion (see e.g. \cite[Lemma 2.1]{khalil_nonlinear_2001}), the set $M$ contains a periodic orbit of $f_0$. 


\aaa{Since $M$ is a one-dimensional connected manifold, the trajectory of $f_0$ starting from a point $\vec{z_0}\in M$ on this periodic orbit can only return to $\vec{z_0}$ by traversing all points in $M.$ Hence, the periodic orbit coincides with $M$. In view of the fact that $\langle \nabla p_{k_0}, f_0\rangle \leq 0$ as established in (\ref{eq:pk0dot<=0}), it follows that $p_{k_0}$ must be equal to some constant $c$ on $M$. Indeed, if we had $p_{k_0}(\vec{z_1})>p_{k_0}(\vec{z_2})$ for some $\vec{z_1}, \vec{z_2}\in M,$ then the trajectory of $f_0$ starting from $\vec{z_2}$ would not visit $\vec{z_1}$ and this is a contradiction.}

\aaa{Note that the constant $c$ must be nonzero or else, by homogeneity, the polynomial $ p_{k_0}$ would be identically zero, contradicting the definition of $k_0$.}
As a consequence, $$p_{k_0} \text{ and } cW^{\frac {k_0}2}$$
 are two \aaa{nonzero} homogeneous functions of \aaa{degree $k_0$} that are equal on $M$. Since $M$ intersects all the lines passing through the origin, \aaa{and since any homogeneous function $u:\mathbb{R}^2\rightarrow\mathbb{R}$ of degree $k_0$ satisfies $u(\lambda x,\lambda y)=\lambda^{k_0}u(x,y)$ for all $\lambda\in\mathbb{R}$ and all $(x,y)\in\mathbb{R}^2$}, we get that 
 $$p_{k_0}(x, y)= cW^{\frac {k_0}2}(x, y) \quad \forall (x, y) \in \mathbb R^2.$$
This implies the following polynomial identity
$$(x^2+y^2)^{k_0}p_{k_0}^2(x, y) = c^2 (x^4+y^4)^{k_0},$$
which gives a contradiction as $(x, y) = (\sqrt{-1}, 1)$ makes only the left-hand side vanish.
\end{proof}

The vector field in \eqref{eq:poly.vector.field.with.rat.coeffs} is a polynomial of degree 7 in two variables. We leave open the problem of determining the minimum degree of a polynomial vector field with rational coefficients for which the statement of Theorem \ref{thm:poly.vector.field.with.rat.coeffs.but.no.poly.lyap} holds. Note also that although the vector field in \eqref{eq:poly.vector.field.with.rat.coeffs} does not admit a polynomial Lyapunov function, it admits a rational one (i.e., a ratio of two polynomials). We leave the question of determining whether LAS polynomial vector fields with rational coefficients admit a local rational Lyapunov function for future research.  We have recently shown in \cite{ahmadi_algebraic_2018} that one cannot hope for a \aaa{\emph{global}} rational Lyapunov function in general. On the other hand, \cite{ahmadi_algebraic_2018} also shows that if the vector field in question is homogeneous, then asymptotic stability implies existence of a rational Lyaponuv function.

\paragraph{Acknowledgments} We are grateful to the anonymous Associate Editor and the anonymous referees whose constructive feedback has improved our presentation considerably.



\bibliographystyle{ieeetr}
\bibliography{citations}

\begin{thebibliography}{10}

\bibitem{khalil_nonlinear_2001}
H.~K. Khalil, {\em Nonlinear {Systems}}.
\newblock Pearson, 3$^{\textrm{rd}}$~ed., 2001.

\bibitem{parrilo_structured_2000}
P.~A. Parrilo, {\em Structured semidefinite programs and semialgebraic geometry
  methods in robustness and optimization}.
\newblock {PhD} {Thesis}, California Institute of Technology, May 2000.

\bibitem{papachristodoulou_construction_2002}
A.~Papachristodoulou and S.~Prajna, ``On the construction of {Lyapunov}
  functions using the sum of squares decomposition,'' in {\em Proceedings of
  the {IEEE} {Conference} on {Decision} and {Control}}, vol.~3, pp.~3482--3487,
  2002.

\bibitem{henrion_positive_2005}
D.~Henrion and A.~Garulli, eds., {\em Positive {polynomials} in {control}}.
\newblock Lecture {Notes} in {Control} and {Information} {Sciences},
  Springer-Verlag, 2005.

\bibitem{jarvis-wloszek_controls_2003}
Z.~Jarvis-Wloszek, R.~Feeley, W.~Tan, K.~Sun, and A.~Packard, ``Some controls
  applications of sum of squares programming,'' in {\em Proceedings of the
  {IEEE} {Conference} on {Decision} and {Control}}, vol.~5, pp.~4676--4681,
  2003.

\bibitem{chesi_lmi_2010}
G.~Chesi, ``{LMI} {techniques} for {optimization} {over} {polynomials} in
  {control}: {a} {survey},'' {\em IEEE Transactions on Automatic Control},
  vol.~55, pp.~2500--2510, 2010.

\bibitem{henrion_guest_2009}
D.~Henrion and G.~Chesi, ``Guest editorial: {special} issue on positive
  polynomials in control,'' {\em IEEE Transactions on Automatic Control},
  vol.~54, no.~5, pp.~935--936, 2009.

\bibitem{zheng2017exploiting}
Y.~Zheng, G.~Fantuzzi, and A.~Papachristodoulou, ``Exploiting sparsity in the
  coefficient matching conditions in sum-of-squares programming using {ADMM},''
  {\em {IEEE} control systems letters}, vol.~1, no.~1, pp.~80--85, 2017.

\bibitem{kamyar2014polynomial}
R.~Kamyar and M.~Peet, ``Polynomial optimization with applications to stability
  analysis and control---alternatives to sum of squares,'' {\em Discrete and
  Continuous Dynamical Systems, Series B}, 2014.

\bibitem{dsosArxiv}
A.~A. Ahmadi and A.~Majumdar, ``{DSOS} and {SDSOS} optimization: more tractable
  alternatives to sum of squares and semidefinite optimization,'' {\em Preprint
  available at arXiv:1706.02586}, 2017.

\bibitem{kamyar2013solving}
R.~Kamyar, M.~M. Peet, and Y.~Peet, ``Solving large-scale robust stability
  problems by exploiting the parallel structure of {P}olya's theorem,'' {\em
  {IEEE} Transactions on Automatic Control}, vol.~58, no.~8, pp.~1931--1947,
  2013.

\bibitem{ben2015stability}
M.~A. Ben~Sassi and S.~Sankaranarayanan, ``Stability and stabilization of
  polynomial dynamical systems using {B}ernstein polynomials,'' in {\em
  Proceedings of the 18th International Conference on Hybrid Systems:
  Computation and Control}, pp.~291--292, ACM, 2015.

\bibitem{peet_exponentially_2009}
M.~M. Peet, ``Exponentially stable nonlinear systems have polynomial {Lyapunov}
  functions on bounded regions,'' {\em IEEE Trans. Automat. Control}, vol.~54,
  no.~5, pp.~979--987, 2009.

\bibitem{peet2012converse}
M.~M. Peet and A.~Papachristodoulou, ``A converse sum of squares {L}yapunov
  result with a degree bound,'' {\em {IEEE} Transactions on Automatic Control},
  vol.~57, no.~9, pp.~2281--2293, 2012.

\bibitem{rational_stability}
T.~Leth, R.~Wisniewski, and C.~Sloth, ``On the existence of polynomial
  {L}yapunov functions for rationally stable vector fields,'' in {\em
  Proceedings of the IEEE Conference on Decision and Control}, pp.~4884--4889,
  2017.

\bibitem{ahmadi_globally_2011}
A.~A. Ahmadi, M.~Krstic, and P.~A. Parrilo, ``A globally asymptotically stable
  polynomial vector field with no polynomial {Lyapunov} function,'' in {\em
  Proceedings of the {IEEE} {Conference} on {Decision} and {Control}},
  pp.~7579--7580, 2011.

\bibitem{bacciotti_liapunov_2006}
A.~Bacciotti and L.~Rosier, {\em Liapunov {Functions} and {Stability} in
  {Control} {Theory}}.
\newblock Springer Science \& Business Media, 2006.

\bibitem{doria_arnolds_1993}
F.~A. Doria and N.~C.~A. da~Costa, ``On {Arnold}'s {Hilbert} symposium
  problems,'' in {\em Computational {Logic} and {Proof} {Theory}}, Lecture
  {Notes} in {Computer} {Science}, pp.~152--158, Springer, 1993.

\bibitem{arnold_problems_1976}
V.~I. Arnold, ``Problems of present day mathematics, {XVII} ({Dynamical}
  systems and differential equations),'' {\em Proc. Symp. Pure Math.}, vol.~28,
  no.~59, 1976.

\bibitem{tarski_decision_1998}
A.~Tarski, ``A {decision} {method} for {elementary} {algebra} and {geometry},''
  in {\em Quantifier {Elimination} and {Cylindrical} {Algebraic}
  {Decomposition}}, Texts and {Monographs} in {Symbolic} {Computation},
  pp.~24--84, Springer, 1998.

\bibitem{seidenberg_new_1954}
A.~Seidenberg, ``A {new} {decision} {method} for {elementary} {algebra},'' {\em
  Annals of Mathematics}, vol.~60, no.~2, pp.~365--374, 1954.

\bibitem{blondel_unsolved_2004}
V.~D. Blondel and A.~Megretski, eds., {\em Unsolved {Problems} in
  {Mathematical} {Systems} and {Control} {Theory}}.
\newblock Princeton University Press, 2004.

\bibitem{jungers_joint_2009}
R.~Jungers, {\em The {Joint} {Spectral} {Radius}: {Theory} and {Applications}}.
\newblock Lecture {Notes} in {Control} and {Information} {Sciences},
  Springer-Verlag, 2009.

\bibitem{lagarias_finiteness_1995}
J.~C. Lagarias and Y.~Wang, ``The finiteness conjecture for the generalized
  spectral radius of a set of matrices,'' {\em Linear Algebra and its
  Applications}, vol.~214, pp.~17--42, 1995.

\bibitem{hare_explicit_2011}
K.~G. Hare, I.~D. Morris, N.~Sidorov, and J.~Theys, ``An explicit
  counterexample to the {Lagarias}–{Wang} finiteness conjecture,'' {\em
  Advances in Mathematics}, vol.~226, no.~6, pp.~4667--4701, 2011.

\bibitem{jungers_finiteness_2008}
R.~M. Jungers and V.~D. Blondel, ``On the finiteness property for rational
  matrices,'' {\em Linear Algebra and its Applications}, vol.~428, no.~10,
  pp.~2283--2295, 2008.

\bibitem{ahmadi_algebraic_2018}
A.~A. Ahmadi and B.~El~Khadir, ``On {algebraic} {proofs} of {stability} for
  {homogeneous} {vector} {fields},'' 2018.
\newblock Preprint available at arXiv:1803.01877.

\end{thebibliography}

\end{document}